\documentclass[12pt, a4paper]{article}
\usepackage{amsmath, amssymb, latexsym}
\renewcommand{\baselinestretch}{1.15}
 \newtheorem{theorem}{Theorem}[section]
\newtheorem{lemma}[theorem]{Lemma}
\newtheorem{proposition}[theorem]{Proposition}

\newtheorem{remark}[theorem]{Remark}

\newtheorem{definition}[theorem]{Definition}
\font\nine=ptmr at 9pt

\newcommand{\R}{\mathbb{R}}

\newcommand{\N}{\mathbb{N}}

\renewcommand{\title}[1]
{\thispagestyle{empty}
\begin{center}
{\Large \bf #1}
\end{center}}

\newcommand{\authors}[1]
{\begin{center}
\renewcommand{\thefootnote}{\fnsymbol{footnote}}
\setcounter{footnote}{3} {\sc #1 }
\end{center}}

\newcommand{\ack}[1]{\footnote{#1}}

\newcommand{\address}[1]
{\vskip 5ex
\renewcommand{\baselinestretch}{1}
\footnotesize \normalsize
 #1 \\
}

\begin{document}

\title{Rearrangement transformations\\ on general measure spaces}

\authors{
Santiago Boza and Javier Soria\ack{Both authors have been
partially supported by Grants MTM2004-02299 and
2005SGR00556.\\{\sl Keywords:} Rearrangement, Lorentz space, symmetrization.\\{\sl  MSC2000:} 46E30, 46B42.} }

\bigskip

 {\narrower\noindent \textbf{Abstract.} \small{For a general set transformation ${\cal R}$ between two measure spaces, we define  the rearrangement of a measurable function by means of the Layer's cake formula. We study some functional properties of the  Lorentz spaces defined in terms of ${\cal R}$, giving a unified approach to the classical rearrangement, Steiner's symmetrization, the multidimensional case, and the discrete setting of trees.} \par}
\bigskip

\section{Introduction}
Given two measure spaces $(X,\Sigma_{X}, \mu)$ and $(Y,\Sigma_{Y},
\nu)$ we consider a general set  transformation ${\cal R}: \Sigma_{X}\rightarrow  \Sigma_{Y}$. 
We denote by $f_{\cal R}^*$  the rearrangement of a
$\mu-$measurable function $f$ with respect the transformation
$\cal R$ by means of the \lq\lq Layer cake formula" (see \cite{LL}):

\begin{equation}\label{lcf}
f_{\cal R}^*(y)=\int_{0}^{\infty} \chi_{{\cal R}(\{x\in X:|f(x)|>t\})}(y) \ dt,
\end{equation}
whenever it defines a $\nu$-measurable function on $Y$.

For $Y=(0,\infty)$ and ${\cal R}$ the transformation defined by
${\cal R}(E)=(0,|E|)$, where $|E|$ denotes the Lebesgue measure of
a set $E\in \Sigma_X$, we have that $f_{\cal R}^*=f^*$, the usual
decreasing rearrangement of a measurable function $f$ defined on
$X$, which we will refer as the classical case (see \cite{BS} for more information).

Formula (\ref{lcf}) has been used recently to define the rearrangement
of functions with respect to some order in very different contexts:
in   \cite{GS1,GS2}
 a new decreasing rearrangement is defined for functions on
homogeneous trees and in \cite{BaPeSo} a multidimensional
rearrangement is considered for functions on $\R^n$.

The work is organized as follows: in Section 2 we develop the main
results concerning  general rearrangements
from a measure theoretical point of view. In Section 3 we
introduce the weighted Lorentz spaces associated to this general
kind of transformations and also we review some functional
properties for these spaces in two known contexts: the
multidimensional rearrangement and the rearrangement on
homogeneous trees, completing the characterization of normability
 already proved in \cite{BaPeSo,GS1}. The theory that we develop allows us to unify and extend these kind of results to two
kind of preserving measures rearrangements that appear very
frequently in applications: Steiner's symmetrization and spherical
rearrangements.

\section{General rearrangement transformations}
In this section, we review how the basic results, well-known in the classical theory, actually imply some a priori assumptions on ${\cal R}$ which, in many cases, turn out to be equivalent statements.
From (\ref{lcf}), we observe that the
rearrangement of a function is a non degenerate transformation;
that is, $f\not\equiv 0$ implies $f_{\cal R}^*\not\equiv 0$, if
there exists $F \in \Sigma_{Y}$ with $\nu(F)>0$ such that a.e.\ 
$y\in F$, there exists $A_y \in (0,\infty)$, with positive Lebesgue
measure, such that $y\in {\bigcap_{t\in A_y}} {\cal
R}\{|f|>t\} $. It is clear from the definitions that having a non
degenerate transformation implies that $\nu({\cal
R}(\emptyset))>0$, or $\nu({\cal R}(E))>0$ if $\mu(E)>0$. The
reverse property is also true if ${\cal R}$ is a monotone
transformation, in the sense that $E\subset F$ implies ${\cal
R}(E)\subset {\cal R}(F)$.

To show that more conditions, like monotonicity, are necessary to have a non
degenerate transformation, let us consider the following
counterexample: assume that $X$ is a subset of $\R^n$ of finite
measure and $Y=(0,\infty)$, with ${\cal R}(E)=(|E|,2|E|)$ (here
$|E|$ denotes the Lebesgue measure of the set $E\subset X$), which is not a monotone rearrangement. An easy application of
the Layer's cake formula shows that $f_{\cal
R}^*(t)=f^*(t/2)-f^*(t)$, where $f^*$ denotes the usual
rearrangement of $f$ with respect to the Lebesgue measure. We deduce
then that any constant function $f$ has $f^*_{\cal R}(y)\equiv 0$,
and so the transformation ${\cal R}$ is degenerate, although in this case, $\nu({\cal R}(E))>0$ if $\mu(E)>0$.
\medskip
\begin{remark}{\rm
\label{simple} A simple non-negative function $f$
can be written as $f(x)={\sum_{k=1}^N b_k
\chi_{F_k}(x)}$, with $(b_k)_{k}>0$ and $(F_k)_{1\leq k\leq N}$ an increasing
sequence of sets. In this case, (\ref{lcf}) gives us that
 \begin{equation}\label{rsf}
 f_{\mathcal{R}}^*(y)=\sum_{k=1}^N b_k \chi_{\mathcal{R}(F_k)}(y),\quad\text {a.e. }y\in Y,
 \end{equation}
provided that the transformation $\mathcal{R}$ satisfies $\nu (
\mathcal{R}(\emptyset))=0.$ In fact, this condition is necessary for (\ref{rsf}) to hold. Thus, from now on, we will always assume that $\nu (\mathcal{R}(\emptyset))=0$, and hence $(\chi_A)^*_{\cal R}=\chi_{{\cal R}(A)}$, a.e., for all $A\in \Sigma_X$.

Also, assuming in addition that ${\cal R}$ is a monotone
transformation on $\Sigma_X$, for
$f(x)={\sum_{j=1}^{N} a_j \chi_{E_j}(x)},$ with
$a_1>a_2>\cdots>a_N>0$ and $E_i\cap E_j=\emptyset$, $i\neq j$,
then
$$
f^*_{\mathcal{R}}(y)=\sum_{j=1}^{N} a_j \chi_{\mathcal{R}(F_j)\setminus
\mathcal{R}(F_{j-1})}(y), \quad\text {a.e. }y\in Y,
$$ 
where
${F_j=\bigcup_{k=1}^{j} E_k}$, and $F_0=\emptyset$.}
 \end{remark}
 
 \begin{lemma}\label{lemmon}
 Suppose $f, g$  are
measurable functions on $X$, and ${\cal R}$ is a monotone transformation. 
\begin{enumerate}
\item [(a)]  If $|g|\leq
|f|$ a.e., then  $ g_{\mathcal R}^*\leq f_{\mathcal R}^*$.
\item[(b)] $ \{f_{\cal R}^*>t\}\subseteqq \mathcal{R}\{|f|>t\}\subseteqq \{f_{\cal R}^* \geq
t\}.$
\item [(c)]  
$(|f|^p)_{\mathcal{R}}^*=(f_{\mathcal{R}}^*)^p,\ 0<p<\infty.$
\end{enumerate}
\end{lemma}
\medskip
\noindent{\bf Proof:} (a) is a consequence of the monotonicity
property on $\mathcal{R}$, because if $|g|\leq |f| $ then $
\chi_{\mathcal{R}(|g|>t)}\leq \chi_{\mathcal{R}(|f|>t)}$ and hence
$ g_{\mathcal R}^*\leq f_{\mathcal R}^*$.
\medskip

To show (b), fix $y\in Y$ such that $f_{\cal R}^*(y)>t$, which is equivalent to  
$$ \int_{0}^{\infty} \chi_{{\cal R}(|f|>s)}(y) \
ds>t. $$ 
Then, $|\{s\in (0,\infty): \ y\in \mathcal{R}(|f|>s)\}|>t$;
i.e., there is $\varepsilon_y>0$ such that 
$$
(0,t+\varepsilon_y)\subset
\{s\in (0,\infty): \ y\in \mathcal{R}(|f|>s)\},
$$ 
and hence $y \in
\mathcal{R}(|f|>t)$. Therefore, we have proved that:
\begin{equation}
\label{1}
 \{f_{\cal R}^*>t\}\subseteqq \mathcal{R}\{|f|>t\}.
 \end{equation}
 On the other hand, since  $y\in \mathcal{R}(|f|>t)$ then $ y\in
\mathcal{R}(|f|>s), \ \forall s \in [0,t]$ and hence $ f_{\cal
R}^*(y) \geq t. $ Therefore,
 \begin{equation}
 \label{2}
\mathcal{R}\{|f|>t\}\subseteqq \{f_{\cal R}^* \geq t\}.
\end{equation}

 To see (c), we use (\ref{1})  and obtain
 \begin{align*}
 (|f|^p)_{\mathcal{R}}^*(x)&=\int_{0}^{\infty} \chi_{{\cal
R}(|f|^p>t)}(x) \ dt=p \int_{0}^{\infty} t^{p-1} \chi_{{\cal
R}(|f|>t)}(x) \ dt \\
& \geq p \int_{0}^{\infty} t^{p-1} \chi_{\{f_{\cal R}^*>t\}}(x) \
dt=p\int_{0}^{f_{\cal R}^*(x)} t^{p-1}  \
dt=(f_{\mathcal{R}}^*)^p(x).
\end{align*}
On the other hand, the inclusion (\ref{2}) establishes the reverse
inequality
\begin{align*}
 (|f|^p)_{\mathcal{R}}^*(x)&=p \int_{0}^{\infty} t^{p-1} \chi_{{\cal
R}(|f|>t)}(x) \ dt \\
& \leq p \int_{0}^{\infty} t^{p-1} \chi_{\{f_{\cal R}^*\geq
t\}}(x) \ dt=p\int_{0}^{f_{\cal R}^*(x)} t^{p-1}  \
dt=(f_{\mathcal{R}}^*)^p(x).
\end{align*}
$\hfill\Box$\bigbreak

\medskip
\begin{definition}
We will say that the transformation ${\cal R}$ satisfies the Fatou
property if  for every increasing sequence of positive 
measurable functions $f_n$ converging to $f$, $\mu$-a.e., we have
that also $(f_n)_{\cal R}^*$ converges increasingly  to $(f)_{\cal
R}^*$, $\nu$-a.e.
\end{definition}
 The following
lemma proves that the Fatou property is equivalent to the fact
that the transformation ${\cal R}$ preserves increasing sequences
of sets. The notation $(A_j)_j\uparrow A$ used below means that $A_j\subset A_{j+1}$, and $A=\cup_j A_j$.

\medskip
\begin{proposition}
\label{FP} Let $(X,\Sigma_X, \mu)$ and $ (Y,\Sigma_Y, \nu)$ be  two
measure spaces and ${\cal R}:\Sigma_X\rightarrow \Sigma_Y$.
Then, the following statements  are equivalent:
\begin{enumerate}
\item [(a)] ${\cal R}$ satisfies the Fatou property.
\item [(b)] For every increasing sequence of sets $(A_j)_j\uparrow A$,
 $(\mathcal{R}(A_j))_j\uparrow \mathcal{R}(A)$. In particular, ${\cal R}$ is monotone.
 \item [(c)]  For $f$ and $f_n$, $n\geq 1$, measurable
functions on $X$,
$$ |f|\leq \liminf_{n\to \infty} |f_n|\Rightarrow f_{\cal R}^*\leq
\liminf_{n\to \infty}(f_n)^*_{\cal R}.$$
\end{enumerate}
\end{proposition}
\medskip

\noindent{\bf Proof:}

First, assume that ${\cal R}$ satisfies the Fatou property. The
condition $(A_j)_j\uparrow A$ is equivalent to
$\chi_{(A_j)}(x)\uparrow \chi_A (x)$ for every $x\in X$ and hence, using (\ref{rsf})
$$(\chi_{A_j})_{\cal R}^*(y)=\chi_{{\cal R}(A_j)}(y)\uparrow \chi_{{\cal R}(A)}(y).$$

To prove that condition (b) implies (c), we define for a fixed
$\lambda\geq 0$,
$$E:=\{x: |f(x)|>\lambda\}, \ \ E_n:=\{x: |f_n(x)|>\lambda\} \
(n=1 ,2,\cdots).$$ Clearly, $E\subset \bigcup_{m\geq 1}
\bigcap_{n\geq m} E_n$ and hence, using (b),
$$
{\cal R}(E)\subseteq {\cal R}\bigg(\bigcup_{m\geq 1} \bigg( \bigcap_{n\geq m} E_n
\bigg) \bigg)= \bigcup_{m\geq 1} {\cal R} \bigg( \bigcap_{n\geq
m} E_n \bigg) \subseteq \bigcup_{m\geq 1}\bigcap_{n\geq m} {\cal
R}(E_n)= \liminf_n {\cal R}(E_n).$$ This inclusion implies that,
for $y\in Y$,
$$
\chi_{{\cal R}(E)}(y) \leq \liminf_{n\to \infty}\chi_{{\cal
R}(E_n)}(y), 
$$ 
and then, using (\ref{lcf}) and
Fatou's lemma
$$f^*_{{\cal R}}\leq \liminf_{n\to
\infty} (f_n)^*_{\cal R}.$$

Finally, it is easy to see that (c) implies that ${\cal R}$ is monotone. Now consider $f_j\uparrow f$. On the one hand, by Lemma~\ref{lemmon}, $(f_j)^*_{\cal R}\le f^*_{\cal R}$ so that  $\limsup_j(f_j)^*_{\cal R}\le f^*_{\cal R}$, and by hypothesis, $f^*_{\cal R}\le \liminf_j(f_j)^*_{\cal R}$, which proves (a).   $\hfill\Box$\bigbreak

\medskip

\begin{theorem} \label{HL} Let ${\cal R}$ be a set transformation between
two measure spaces $(X,\Sigma_X, \mu)$ and $ (Y,\Sigma_Y, \nu)$.
Assume that ${\cal R}$ satisfies the Fatou property.  Then, the
following are equivalent conditions:
\begin{enumerate}
\item [(a)] $\mu(A\cap B)\leq \nu({\cal R}(A)\cap {\cal R}(B))$, for every $A,B\in\Sigma_X$.
\item [(b)] $\int_A f \ d\mu \leq \int_{{\cal R}(A)} f_{{\cal R}}^*  \
d\nu$, for every  non-negative measurable function $f$ on $X$, and $A\in\Sigma_X$.
\item [(c)] $\int_X f g \ d\mu \leq \int_Y f_{{\cal R}}^* g_{{\cal R}}^* \  
d\nu$, for every  non-negative measurable functions  $f$ and $g$ on $X$.
\end{enumerate}
\end{theorem}
\medskip

\noindent{\bf Proof:} Let us assume (a). Using the Fatou property,
it is enough to prove (b) just for a simple function of the form
$$
f(x)=\sum_{j=1}^N b_j \chi_{E_j}(x), \ \ \mbox{with $(b_j)_j>0$ and $E_j$ an increasing sequence of sets},
$$
since we can always find a
sequence $(s_k)_k$ of simple functions such that
$0\leq (s_1)_{\mathcal R}^*\leq \cdots\leq (s_k)_{\mathcal{R}}^*
\leq f_{\mathcal{R}}^*$
and $(s_k)_{\mathcal{R}}^*(y)\to f_{\mathcal{R}}^*(y)$, as $k\to
\infty$, a.e.\ $y\in Y$.

Then,
\begin{align*}
\int_A g \ d\mu &=\sum_{j=1}^N b_j \mu(A\cap E_j) \leq \sum_{j=1}^N b_j \nu({\cal R}(A)\cap {\cal R}(E_j))  \\
&=\sum_{j=1}^N b_j \int_{{\cal R}(A)} \chi_{{\cal R}(E_j)}(y) \ \
d\nu(y) =\int_{{\cal R}(A)} g_{\cal R}^* \ d\nu.
\end{align*}
To prove (c) assuming (b), we can also suppose $f(x)=\sum_{j=1}^{N}
a_j \chi_{E_j}(x)$, $a_j>0$ and $E_j$ an increasing sequence of
sets. Then, by (b),
\begin{align*}
\int_X fg \ d\mu &=\sum_{j=1}^N a_j \int_{E_j} g \ d\mu \leq \sum_{j=1}^N a_j \int_{{\cal R}(E_j)} g_{\cal R}^* \ d\nu  \\
&=\int_Y \sum_{j=1}^N a_j \chi_{{\cal R}(E_j)}(y) g_{\cal R}^*(y)
\ d\nu= \int_{Y} f_{{\cal R}}^*(y) g_{{\cal R}}^*(y) \ d\nu.
\end{align*}

Finally if we take $f=\chi_A$ and $g=\chi_B$ in condition (c) we
obtain (a). $\hfill\Box$\bigbreak

\medskip
\begin{definition}
We will say that ${\cal R}$ is a measure preserving  transformation from $\Sigma_X$  into
$\Sigma_Y$, if $\mu(E)=\nu({\mathcal R(E)})$, for every $E \in
\Sigma_X$.
\end{definition}

\begin{proposition}
\label{preservtransf} Let us suppose that ${\mathcal R}$ is a
monotone transformation. Then, the following statements are equivalent:
\begin{enumerate}
\item [(a)] ${\mathcal R}$ is a measure preserving transformation.

\item [(b)] If  $s(x)={\sum_{j=1}^{N} a_j \chi_{E_j}(x)},$ with
$a_1>a_2>\cdots>a_N>0$, $E_i\cap E_j=\emptyset$, $i\neq j$ and  $p>0$, then
$$ 
\int_X s(x)^p \ d\mu(x)= \int_Y
(s_{\mathcal{R}}^*)^p(y) \ d\nu(y).
$$
\end{enumerate}
\end{proposition}

\medskip
\noindent{\bf Proof:} If ${\cal R}$ is a monotone measure preserving transformation,  and $E\cap F=\emptyset$,  with  $\mu(E)<\infty$, then 
$$
\nu(\mathcal{R}(E\cup F)\setminus \mathcal{R}(E)) =\nu
(\mathcal{R}(E\cup F))- \nu(\mathcal{R}(E))=\mu(E\cup
F)-\mu(E)=\mu(F). 
$$
Thus, (b) follows by Remark~\ref{simple}, since$$
(s_{\mathcal{R}}^*)^p(y)=\sum_{j=1}^N a_j^p \chi_{\mathcal{R}(F_j)\setminus
\mathcal{R}(F_{j-1})}(y),
$$
 with
${F_j=\bigcup_{k=1}^{j} E_k}$, and
$F_0=\emptyset$. (a) follows from (b) by taking $s=\chi_A$.$\hfill\Box$\bigbreak

\section{Lorentz spaces and symmetrization}
In this section we prove some properties of a new type of Lorentz
spaces, defined using the general transformations ${\cal R}$. Let
$v$ be a weight on $Y$ (i.e., $v\in L^1_{\hbox{\nine loc}}(Y,d\nu)$, $v\ge 0$ and satisfies the following non-cancellation property: if $\mu(A)>0$, then $\int_{{\cal R}(A)}v(y)\,d\nu(y)>0$), and $0<p<\infty$. We will say that a
$\mu$-measurable function on $X$ belongs to the Lorentz space
$\Lambda^p_{\cal R}(v)$, provided $\|f\|_{\Lambda^p_{\cal R}(v)}$,
defined by
\begin{equation} \label{lornorm}
 \|f\|_{\Lambda^p_{\cal R}(v)}:=\left(\int_Y
(f_{\cal R}^*(y))^p \ v(y) \ d\nu(y) \right)^{1/p},
\end{equation}
 is finite. The case $X=\mathbb R^n$, $Y=\mathbb R^+$, ${\cal R}(E)=(0,|E|)$, and $v(y)=y^{p/q-1}$ gives the classical Lorentz space: $\Lambda^p_{\cal R}(v)=L^{q,p}(\mathbb R^n)$.

The question whether the functional defined in (\ref{lornorm})
is a norm  was answered by Lorentz in the euclidean case (see \cite{Lo} for a proof and \cite{CS3,Saw2} for related questions). Also, M.J.\ Carro and J.\ Soria (\cite{CS}) characterized
the weights $v$ such that it becomes a quasi-norm, if $X$ is no
atomic. Later, in \cite{CRS}, the quasi-normability was completed for all $X$.
The analogous characterization was established in \cite{BaPeSo}
for the multidimensional rearrangement and in \cite{GS1} for the
case of homogeneous trees. In this section we give partial answers
to this question in the context of a general  transformation
${\cal R}$, satisfying the Fatou property, between $\sigma$-finite measure spaces
$X$ and $Y$ (from now on, we will always assume these two conditions).

We adopt the notation $V(E)=\int_{E} v(y) \ d\nu(y)$, for every
measurable set $E\subset Y$ and every weight $v$ in $Y$. Then, the
functional (\ref{lornorm}) has the following description:

\begin{lemma}
\label{labelset} Let $0<p<\infty$. Then, for all $f\in
\Lambda^p_{\cal R}(v)$, we have
\begin{equation}
\|f\|_{\Lambda^p_{\cal R}(v)}=\left(\int_0^{\infty} p
\lambda^{p-1} V({\cal R}(|f|>\lambda)) \  d\lambda \right)^{1/p}.
\end{equation}
\end{lemma}
\medskip
\noindent{\bf Proof:} Using Lemma~\ref{lemmon} (c) we have:
$$\|f\|_{\Lambda^p_{\cal R}(v)}=\left(\int_Y (|f|^p)_{\cal R}^*(y)^p \
v(y) \ d\nu(y) \right)^{1/p}.$$ Then,  by (\ref{lcf}) and Fubini's Theorem,
\begin{eqnarray*}
 \|f\|_{\Lambda^p_{{\cal
R}(v)}} &=& \left(\int_Y \left(\int_0^{\infty} \chi_{{\cal
R}(|f|^p>\lambda)}(y) \ d\lambda\right) \ v(y) \ d\nu(y)
\right)^{1/p} \\
&=&\left(\int_Y \left(\int_0^{\infty} p \xi^{p-1} \chi_{{\cal
R}(|f|>\xi)}(y) \ d\xi \right) \ v(y) \ d\nu(y) \right)^{1/p} \\
&=&\left(\int_0^{\infty} p \xi^{p-1}\left(\int_{ {\cal
R}(|f|>\xi)} v(y) \ d\nu(y) \right) \ d\xi\right)^{1/p}.
\end{eqnarray*}
$\hfill\Box$\bigbreak

Our first result gives a characterization of the quasi-normability of the
functional defined in (\ref{lornorm}).

\begin{theorem}
The functional $\|\cdot\|_{\Lambda^p_{\cal R}(v)}$ is a quasi-norm
if and only if there exists a constant $C>0$ such that
\begin{equation}
\label{quasi}
 V({\cal R}(A\cup B))\leq C (V({\cal R}(A))+V({\cal
R}(B))),
\end{equation}
 for all sets $A,B\in\Sigma_X$.
\end{theorem}
\medskip \noindent{\bf Proof:} Assume first (\ref{quasi}): by
Lemma \ref{labelset}, if $\|f\|_{\Lambda^p_{\cal R}(v)}=0$, then
$$V({\cal R}\{|f|>\lambda\})=0,$$
for all $\lambda>0$, and by hypothesis
$\mu(\{|f|>\lambda\})=0$, for all $\lambda$; that is $f\equiv
0$. Also by Lemma \ref{labelset}, the hypothesis and   the monotonicity in ${\cal R}$, we have:
\begin{eqnarray*}
\|f+g\|^p_{\Lambda^p_{\cal R}(v)}&=&\int_0^{\infty} p
\lambda^{p-1} V({\cal R}(|f+g|>\lambda)) \  d\lambda \\
&\leq & \int_0^{\infty} p \lambda^{p-1} V({\cal
R}(\{|f|>\lambda/2\}\cup \{|g|>\lambda/2\})) \  d\lambda \\
&\leq& C \left(\int_0^{\infty} p \lambda^{p-1} V({\cal
R}(|f|>\lambda/2)) \  d\lambda+\int_0^{\infty} p \lambda^{p-1}
V({\cal R}(|g|>\lambda/2)) \  d\lambda \right)\\
&=& 2C \left(\int_0^{\infty} p \lambda^{p-1} V({\cal
R}(|f|>\lambda)) \  d\lambda+\int_0^{\infty} p \lambda^{p-1}
V({\cal R}(|g|>\lambda)) \  d\lambda \right)\\
&=& 2C (\|f\|^p_{\Lambda^p_{\cal R}(v)}+\|g\|^p_{\Lambda^p_{\cal
R}(v)})\le C_p(\|f\|_{\Lambda^p_{\cal R}(v)}+\|g\|_{\Lambda^p_{\cal
R}(v)})^p.
\end{eqnarray*}
Conversely, suppose that the functional is a quasi-norm and take
$A$ and $B$. Then,
\begin{eqnarray*}
V({\cal R}(A\cup
B))^{1/p}&=&\|\chi_{A\cup B}\|_{\Lambda^p_{\cal R}(v)} \leq
C(\|\chi_A\|_{\Lambda^p_{\cal R}(v)}+ \|\chi_B\|_{\Lambda^p_{\cal
R}(v)})\\
&=& C (V({\cal R}(A))^{1/p}+V({\cal R}(B))^{1/p})\le C (V({\cal R}(A))+V({\cal R}(B)))^{1/p}. 
\end{eqnarray*}
$\hfill\Box$\bigbreak

Concerning the normability of $\Lambda^p_{\cal R}(v)$, we can
establish the following partial results:
\begin{theorem}\label{tconcavity}
 Let $1\leq
p<\infty$, and $v$ be a weight on $Y$. If
$\|\cdot\|_{\Lambda^p_{\cal R}(v)}$ is a norm then, for all $A,B\in\Sigma_X$,
\begin{equation}\label{concavity}
V({\cal R}(A\cup B))+V({\cal R}(A\cap B))\leq  V({\cal
R}(A))+V({\cal R}(B)).
\end{equation}
\end{theorem}
\medskip
\noindent{\bf Proof:} If $\|\cdot\|_{\Lambda^p_{\cal R}(v)}$ is a
norm, take $A$,$B \subset X$, $\delta>0$ and define the functions
$$
f(x)=(1+\delta)\chi_{A}(x)+\chi_{(A\cup B)\setminus A}(x)
$$
and
$$
g(x)=(1+\delta)\chi_{A}(x)+\chi_{(A\cup B)\setminus B}(x).
$$
Then,
\begin{eqnarray*}
f_{\cal{R}}^*(y)&=&(1+\delta)\chi_{{\cal R}(A)}(y)+\chi_{{\cal R}(A\cup B)\setminus {\cal R}(A)
}(y),\\
g_{\cal{R}}^*(y)&=&(1+\delta)\chi_{{\cal
R}(B)}(y)+\chi_{{\cal R}(A\cup B)\setminus {\cal R}(B) }(y),\\
(f+g)_{\cal{R}}^*(y)&=&(2+2\delta)\chi_{{\cal R}(A\cap B)}(y)+(2+\delta) \chi_{{\cal R}(A\cup B)\setminus {\cal R}(A\cap B)
}(y).
\end{eqnarray*}
 The triangle inequality and the fact that $1/p\leq 1$ imply
\begin{eqnarray*}
 \|f+g\|_{\Lambda^p_{{\cal R}(v)}}&=&\left((2+2\delta)^p V({\cal
R}(A\cap B))+(2+\delta)^p V({\cal R}(A\cup B)\setminus {\cal
R}(A\cap B))\right)^{1/p}\\
& \leq & \|f\|_{\Lambda^p_{{\cal R}(v)}}+\|g\|_{\Lambda^p_{{\cal
R}(v)}}=\left( (1+\delta)^p V({\cal R}(A))+V({\cal R}(A\cup
B)\setminus {\cal R}(A)) \right)^{1/p}\\
&\qquad+&\left( (1+\delta)^p V({\cal R}(B))+V({\cal R}(A\cup B)\setminus
{\cal R}(B)) \right)^{1/p}\\
&\leq & 2^{1-1/p} ((1+\delta)^p V({\cal R}(A))+V({\cal
R}(A\cup B)\setminus {\cal R}(A)) \\
&\qquad+& (1+\delta)^p V({\cal R}(B))+V({\cal R}(A\cup B)\setminus {\cal
R}(B)))^{1/p}.
\end{eqnarray*}

Collecting terms, dividing both sides by $2^{p-1}((1+\delta)^p-1)$
and letting $\delta \to 0$, we finally obtain
$$V({\cal R}(A\cup B))+V({\cal R}(A\cap B))\leq  V({\cal
R}(A))+V({\cal R}(B)).$$ $\hfill\Box$ \bigbreak

Condition (\ref{concavity}), in the classical case, implies that $V$ is a concave function, and we will refer to it as the {\sl Concavity Condition.} A sufficient condition in a general setting to ensure that
the functional $\|\cdot\|_{\Lambda^p_{\cal {R}}(v)}$ defines a
norm is the following {\sl Saturation
Property:}

\begin{theorem}
\label{saturation} Let $1\leq p<\infty$, and $v$ be a weight on $Y$
such that $v$ coincides with $h_{\cal R}^*$ for some $h$ defined
on $X$. Then, $\|\cdot\|_{\Lambda^p_{\cal R}(v)}$ is a norm if for
all measurable functions $f$ in $X$, the equality
\begin{equation}\label{sp}
\sup_{\{h:h_{\cal R}^*=v\}}\int_X |f(x) \ h(x)| \
d\mu(x)=\int_{Y} f_{\cal R}^*(y) \ v(y) \ d\nu(y) 
\end{equation} 
holds.
\end{theorem}
\medskip
\noindent{\bf Proof:} We apply  Lemma~\ref{lemmon} (c) and
the hypothesis:
\begin{align*}
\|f+g\|_{\Lambda^p_{\cal R}(v)}&=\left(\int_Y (f+g)_{\cal
R}^{*p}(y)\ v(y) \ d\nu(y) \right)^{1/p}=\left(\int_Y
(|f+g|^p)_{\cal R}^{*}(y)\ v(y) \ d\nu(y) \right)^{1/p}\\
&=\sup_{\{h:h_{\cal R}^*=v\}} \left( \int_X |f(x)+g(x)|^p \ h(x) \
d\mu(x) \right)^{1/p} \\
&\leq \sup_{\{h:h_{\cal R}^*=v\}} \left(
\int_X |f(x)|^p \ h(x) \ d\mu(x) \right)^{1/p}+\sup_{\{h:h_{\cal
R}^*=v\}} \left( \int_X
|g(x)|^p \ h(x) \ d\mu(x) \right)^{1/p}\\
&=\left(\int_Y (f_{\cal R})^{*p}(y)\ v(y) \ d\nu(y) \right)^{1/p}
+\left(\int_Y (g_{\cal R})^{*p}(y)\ v(y) \ d\nu(y)
\right)^{1/p}\\
&=\|f\|_{\Lambda^p_{\cal R}(v)}+\|g\|_{\Lambda^p_{\cal R}(v)}.
\end{align*}
$\hfill\Box$ \bigbreak

\begin{remark}{\rm 
We observe that in order for $\|\cdot\|_{\Lambda^p_{\cal R}(v)}$ to be a norm is not enough 
that the weight $v$ be the rearrangement of some function $h$
defined on $X$. The conditions are, in general, more restrictive:
see \cite{GS1} in the case of trees or \cite{BaPeSo} in the
multidimensional setting. In the next section we will deal with
these examples.}
\end{remark}

Even though normability can fail,  completeness of ${\Lambda^p_{\cal R}(v)}$ always holds:

\begin{proposition}
Assume that $v$ is a weight on $Y$, such that the Lorentz
space $\Lambda:=\Lambda^p_{\cal {R}}(v)$ is continuously embedded
in the space $L_{\hbox{\nine loc}}^1(X)$, and $\|\cdot\|_{\Lambda}$ is a quasi-norm. If $(f_n)$ is a Cauchy sequence
in $\Lambda$ then, there exists a measurable function $f\in
\Lambda$ such that $\lim_{n\to\infty} \|f-f_n\|_{\Lambda}=0$.
\end{proposition}
\medskip
 \noindent{\bf
Proof:} Since $\|\cdot\|$ is a quasi-norm and $(f_n)$ is Cauchy,
there exists a constant $C>0$ such that $\|f_n||_{\Lambda}\leq
C<\infty$, for all $n\in \N$.

Since $X$ is $\sigma$-finite, let us write
$X={\bigcup_{k\geq 1} A_k}$, with $\mu(A_k)<\infty$
and $A_k$ an increasing sequence of sets.

It is clear that $f_n
\chi_{A_k}$ is a Cauchy sequence in $L^1(A_k)$ and hence the sequence $f_n
\chi_{A_k}$ converges to a function $g_k$ in $L^1(A_k)$, for each
$k$. Let us define $f:=g_k$ in $A_k$, which is well-defined by the monotonicity of $A_k$. We have to prove that $f_n
\to f $ in $\Lambda$. By standard arguments, we can find a subsequence $f_{j_k}\to f $ a.e. $x\in X$. Then,  by  Proposition \ref{FP} (c) and Fatou's lemma, we have that
$f\in \Lambda$, and 
\begin{eqnarray*}
\int_Y (f_{\cal R}^*)^p(y) \ v(y) \ d\nu(y) &\leq& \int_{Y}
\liminf_k (f_{j_{k}})_{\cal R}^{*p} \ v(y) \ d\nu(y) \\
&\leq & \liminf_k \int_Y (f_{j_k})_{\cal R}^{*p} \ v(y) \
d\nu(y)=\liminf_k \|f_{j_k}\|_{\Lambda}^p \leq C^p.
\end{eqnarray*}

Using Fatou's lemma again and the fact that $(f_k)_k$ is a Cauchy
sequence, we finally get
$$\|f-f_n\|_{\Lambda}\leq
C(\|f-f_{j_k}\|_{\Lambda}+\|f_{j_k}-f_n\|_{\Lambda})\to 0,$$ as
$n,k\to \infty$. $\hfill\Box$\bigbreak

\begin{definition}
A weight $v$ defined on the
space $Y$ is called ${\cal R}$-admissible if for  every $A \in \Sigma_X$ and
all $0<\varepsilon < \int_{{\cal R}(A)} v(y) \ d\nu(y)$, there exists
${\cal R}(B)\subset {\cal R}(A), $ such that $\int_{{\cal R}(B)}
v(y) \ d\nu(y)=\varepsilon$.
\end{definition}

Now, we can show the following necessary condition on $p$ for
which the functional $\|\cdot\|_{\Lambda^p_{\cal R}(v)}$ defines a
norm:

\begin{theorem}
\label{admissible}
  Let   $v$ be an ${\cal
R}$-admissible weight and  $0<p<\infty$. If $\Lambda^p_{{\cal R}}(v)$ is a Banach
space, then $p\geq 1$.
\end{theorem}

\medskip
\noindent{\bf Proof:} Since $\Lambda^p_{{\cal R}}(v)$ is a Banach
space, there exists $\|\cdot \|$, a norm on $\Lambda^p_{{\cal
R}}(v)$, such that $\|f\|_{\Lambda^p_{{\cal R}}(v)}\simeq \|f\|$.
Hence,
$$\bigg\|\sum_{k=1}^{N} f_k \bigg\|_{\Lambda^p_{{\cal R}}(v)}\leq C
\sum_{k=1}^N \|f_k\|\leq \tilde{C} \sum_{k=1}^N
\|f_k\|_{\Lambda^p_{{\cal R}}(v) },$$ for all $N\in \N$. Suppose
$0<p<1$. Due to the hypothesis assumed on $v$, we can take a
decreasing sequence of subsets
$$A_{k+1}\subset A_{k} \cdots \subset X,$$
such that $\int_{{\cal R}(A_k)}v(y) \ d\nu(y)=2^{-kp}$. If
$f_k=2^k \chi_{A_k}$, then $\|f_k\|_{\Lambda^p_{{\cal R}}(v)}=1.$
But for a fixed $N$, we have
$$\frac{1}{N}\bigg\|\sum_{k=1}^{N} f_k \bigg\|_{\Lambda^p_{{\cal R}}(v)}\leq
\tilde{C}<\infty.$$ On the other hand, since by Remark~\ref{simple}
$(\sum_{k=1}^N 2^k \chi_{A_k})_{{\cal R}}^*=\sum_{k=1}^N 2^k
\chi_{{\cal R}(A_k)}$ and ${\cal R}(A_{k+1})\subset {\cal
R}(A_k)\subset \cdots \subset Y$, we have (taking ${\cal
R}(A_{N+1})=\emptyset$),
\begin{eqnarray*}
& &\frac{1}{N}\bigg\|\sum_{k=1}^{N} f_k \bigg\|_{\Lambda^p_{{\cal
R}}(v)}=\frac{1}{N}\bigg\|\sum_{k=1}^{N} 2^k \ \chi_{A_k}
\bigg\|_{\Lambda^p_{{\cal R}}(v)}\\
&=&\frac{1}{N}\left(\int_{Y} \left(\sum_{k=1}^N 2^k \chi_{{\cal
R}(A_k)}\right)^p(y) \ v(y) \ d\nu(y) \right)^{1/p}\\ &=&
\frac{1}{N}\left(\int_{Y} \left(\sum_{k=1}^N  \left(\sum_{j=1}^k
2^j\right)^p \chi_{{\cal R}(A_k)\setminus {\cal
R}(A_{k+1})}(y)\right)
 \ v(y) \ d\nu(y) \right)^{1/p}\\
 &=&\frac{1}{N}\left(\sum_{k=1}^N  \left(\sum_{j=1}^k
 2^j\right)^p \left(\int_{{\cal R}(A_k)}v(y) \ d\nu(y)-\int_{{\cal
R}(A_{k+1})}v(y) \ d\nu(y)\right) \right)^{1/p}\\  & \geq &
\frac{C}{N} \left(\sum_{k=1}^N(1-2^{-k})^p\right)^{1/p} \geq
\frac{C}{N}\left(\sum_{k=1}^N
2^{-p}\right)^{1/p}=C\frac{N^{1/p}}{N} \to \infty, \ \mbox{as
$N\to \infty$,}
\end{eqnarray*}
which is a contradiction. Hence $p\geq 1$.$\hfill\Box$ \bigbreak

\begin{remark}{\rm
We observe that in Theorem~\ref{admissible} the hypothesis assumed
on $Y$ is not compatible with the fact that $Y$ is a completely atomic measure space. In the case $0<p<1$, if $Y$ is completely
atomic, we observe that the functional
$\|\cdot\|_{\Lambda^p_{{\cal R}}(v)}$ is a norm if and only if
supp $v$ is contained in some atom ${\cal R}(A)$ such that, for
every measurable set ${\cal R}(B)$ in $Y$, ${\cal R}(A)\subset
{\cal R}(B).$ Observe that this is the case of the discrete
setting (see \cite{GS1} for a proof in the context of homogeneous
trees).}
\end{remark}

The classical Lorentz spaces are generalizations of the
Lebesgue spaces, since $\Lambda_{X}^p(1)=L^p(X)$. The next proposition
shows that for a general transformation  ${\cal R}$, the
corresponding Lorentz space also satisfies this property provided
that ${\cal R}$ is a measure preserving transformation.

\begin{proposition}
Let $0<p<\infty$. Then, ${\cal R}$ is a measure preserving
transformation if and only if $\Lambda^p_{{\cal R}}(1)=L^p(X)$, with equality of norms.
\end{proposition}

\medskip
\noindent{\bf Proof:} If ${\cal R}$ is a measure preserving
transformation, by Fubini's theorem and Lemma~\ref{lemmon}, we have:
\begin{eqnarray*}
\|f\|_{L^p(X)}^p&=& \int_{X}|f(x)|^p\ d\mu(x)=\int_0^{\infty}
\int_{\{|f|^p>t\}} \ d\mu(x) \ dt=\int_0^{\infty} \int_{{\cal
R}(|f|^p>t)} \ d\nu(y) \ dt\\&=&\int_Y \int_0^{\infty} \chi_{{\cal
R}(|f|^p>t)}(y) \ dt \ d\nu(y) =\int_Y (|f|^p)^*_{{\cal R}}(y) \
d\nu(y)\\ &=&\int_Y (|f|_{\cal R}^*)^p(y)\ d\nu(y)
=\|f\|_{\Lambda^p_{{\cal R}}(1)}^p.
\end{eqnarray*}
The converse follows by taking $f=\chi_A$. $\hfill\Box$ \bigbreak

In the general context of a monotone transformation ${\cal R}$
between measure spaces, Theorems~\ref{tconcavity} and
\ref{saturation} give two conditions (one necessary and the other sufficient)  to ensure that
the functional given by (\ref{lornorm}) defines a norm. Both
conditions are known as the concavity condition and the saturation
property, respectively, and are equivalent in the   classical setting. Moreover, 
they are also  equivalent to the fact that the weight $v$ must be
decreasing (see \cite{Lo}).

In the case of the two-dimensional rearrangement it has been
proved that (\ref{lornorm}) is a norm if and only if the concavity
condition holds and the weight $ v$ defined on $\R^2_+$ is a
decreasing function that only depends on one variable (see
\cite[Theorem 3.7]{BaPeSo}). On the other hand, in the case of
rearrangement defined on homogeneous trees it has been shown (see
 \cite[Theorem 4.9]{GS1}) that the saturation property holds for
linear decreasing weights (see \cite{GS2} for the definition) and
both conditions are equivalent to the fact that (\ref{lornorm})
defines a norm.

We can briefly resume these conditions in the following list:
\begin{align*}
\text{{\bf (Norm):} (\ref{lornorm}) defines a norm. \ \ \ }
&\text{{\bf(CC):}  Concavity Condition  (\ref{concavity}).}\\
\text{{\bf(SP):} Saturation Property  (\ref{sp}). }
&\text{{\bf(MP):} Monotonicity properties on the weight.}
\end{align*}

Then, in the classical setting:
$$
\text{{\bf (Norm)}}\iff\text{{\bf (CC)}}\iff\text{{\bf (SP)}}\iff\text{{\bf (MP)}},
$$
 in the multidimensional  setting:
$$
\text{{\bf (Norm)}}\iff\text{{\bf (CC)}}\iff\text{{\bf (MP)}}\Longleftarrow\text{{\bf (SP)}},
$$
and in the case of trees:
$$
\text{{\bf (Norm)}}\iff\text{{\bf (SP)}}\iff\text{{\bf (MP)}}\implies\text{{\bf (CC)}}.
$$
We will now complete the missing results in the above list, and extend the equivalences to two more rearrangements (spherical and Steiner's symmetrization).

In the case of the multidimensional rearrangement, to simplify the notation, we will restrict ourselves to the
two-dimensional case. We can establish the following saturation
property which completes the characterization of normability of
Lorentz spaces in this context (see \cite{BaPeSo2,BaPeSo}). Also, it is proved in
\cite{BaPeSo}  that given a function $f(x,y)$ defined
on $\R^2 $, its two dimensional rearrangement, $f_2^*(s,t)$,
$s,t>0$, can be understood as an iterative procedure with respect
to the usual rearrangement in each variable. More precisely,
$f_2^*(s,t)=(f_y^*(\cdot,t))_x^*(s)$. That is, first we rearrange
with respect to $y$ and after with respect to $x$. In this case {\bf (MP)} is given by the fact that the weight $v(s,t)=v(t)$, where $v$ is a decreasing function.

\begin{proposition}
For any measurable function in $\R^2$,
$$\sup_{h_2^*=v} \int_{\R^2} f(x,y) \ h(x,y) \ dx \
dy=\int_{\R_{+}^2} f_2^*(s,t) \ \ v(t) \ ds \ dt,$$ where $v(t)$
is a decreasing function with respect to the variable $t \in \R^+$.
\end{proposition}
\medskip
\noindent{\bf Proof:} Applying Hardy-Littlewood inequality with
respect the one dimensional decreasing rearrangement we have that,
\begin{eqnarray*}
\int_{\R^2} f(x,y) h(x,y) \ dy \ dx&=& \int_{\R} \int_{\R} f_x(y)
\ h(x,y) \ dy \ dx
\leq \int_{\R} \int_0^{\infty} f_y^*(x,t) \ h_y^*(x,t) dt \ dx \\
&\leq&\int_0^{\infty} \int_0^{\infty} (f_y^*(\cdot,t))_x^*(s) \
(h_y^*(\cdot,t))_x^*(s) \ dt\  ds\\
&=&\int_0^{\infty} \int_0^{\infty} f_2^*(s,t) \ v(t) \ dt\ ds.
\end{eqnarray*}

To prove the converse, we use that for $u$ a decreasing
function (see \cite{BS}),
\begin{equation}
\label{satur}
 \sup_{\sigma} \int_{\R^n} |f(x)| \ u(\sigma(x)) \ dx
=\int_0^{\infty} f^*(t) u(t) \ dt,
\end{equation}
 where the supremum is taken
over all measure preserving transformations
$\sigma:\R^n\longrightarrow \R^{+}$.

 Let us show that
\begin{equation}
\label{satur3}
 \int_0^{\infty}\int_0^{\infty} (f_y^*(\cdot, t))_x^*(s) \
v(t) \ dt \ ds \leq \sup_{\sigma} \int_{\R^2} f(x,y) \
v(\sigma_x(y)) \ dy \ dx.
 \end{equation}

For a given $\varepsilon>0$ and $x\in \R$, using (\ref{satur}), there
exists $\sigma_x:\R\longrightarrow \R^{+}$ such that
$$\frac{1}{1+\varepsilon} \int_{\R} f_y^*(x,t)\ v(t) \ dt
\leq \int_{\R} f(x,y)\ v(\sigma_x(y)) \ dy.$$

We integrate over $x\in \R$ and obtain
\begin{eqnarray*}
 \frac{1}{1+\varepsilon} \int_{\R} \left( \int_{0}^{\infty}
f_y^*(x,t) \ v(t) \ dt \right)\ dx &=&\frac{1}{1+\varepsilon}
   \int_0^{\infty} \int_0^{\infty} (f_y^*(\cdot, t))_x^*(s) \
v(t) \ dt \ ds                                              \\
&\leq& \sup_{\sigma} \int_{\R^2} f(x,y) \ v(\sigma_x(y)) \ dy \
dx,
\end{eqnarray*}
which gives us 
(\ref{satur3}). $\hfill\Box$ \bigbreak

In the case of a tree, we can complete the set of
equivalences by showing that also the concavity condition for a
function $v$, defined on the tree, implies that $v$ is a linear
decreasing weight, which is {\bf (MP)} (see \cite{GS1, GS2}):

\begin{definition}
For two given disjoint sets $A$ and $B$ in the boundary $\partial
T$ of a homogeneous tree $T$, we write $A\le B$, if
$\alpha\le\beta$ for all $\alpha\in A$ and all $\beta\in B$. Then, given two
vertices $x$ and $y$ in $T$, we define
$$x\unlhd y$$
if and only if
$$x\leq y\ \ \mbox{or} \ \ I(x)\geq I(y),$$
where $I(x)$ is the set of all geodesics passing through $x$. 
We say that the function $f$ is linearly decreasing if $f(x)\geq
f(y)$ whenever $x\unlhd y$.
\end{definition}

\begin{proposition}
Let $v$ be a weight in $T$. If $v$ satisfies the concavity
condition {\bf (CC)}, then $v$ is linearly decreasing.
\end{proposition}
\medskip
\noindent{\bf Proof:} Let us consider two vertices $x\lhd y$. It is
enough to consider the case $ I(x)\geq
I(y).$

Set $A=[o,x]$ and $B=[1,y]\cup [1,x]$.
If we denote by $x-1$   the vertex in the geodesic $[o,x]$ with distance to $x$ equal to 1 then,  $A=A^*$, $B^*=[o,y-1]\cup [1,x]$, $A\cup B=(A\cup
B)^*=[o,y]\cup[1,x]$ and $(A\cap B)^*=[o,x-1]$.

Applying {\bf (CC)} for these sets  
we easily obtain that $v(y)\leq v(x)$; that is, $v$  is linearly
decreasing.$\hfill\Box$ \bigbreak

We will now consider two more well-known rearrangements. Let $A$ be a
measurable set in $\R^n$, $n\geq 2$. The spherical symmetrization
of a set $A$ is ${\cal R}(A)=A^*=B(0,(\sigma_n^{-1}|A|)^{1/n})$,  where
$\sigma_n$ is the volume of the $n$-dimensional ball (see \cite{AL} for further information). To define
Steiner symmetrization (see \cite{LL,B}) of order $k
\geq 1$, we write points in $x\in \R^{n}$ as pairs $x=(\bar{x},y)$
with $\bar{x}\in \R^{n-k}$ and $y\in \R^{k} $. The Steiner
symmetrization of order $k$ of $A$ is ${\cal R}(A)={\cal S}_k(A)$,  the set
whose $k$-dimensional cross sections parallel to the hyperplane
$\bar{x}=0$ are balls with measure equal to the corresponding
cross sections of $A$. This symmetrization method shows up in applications to PDE's, like the isoperimetric inequality (see \cite{LL,Ka}).

For $f:\R^n \longrightarrow \R$ a measurable function, we define
the spherical symmetrization $f_{Sp}^*$ and the Steiner
symmetrization $({\cal S}_k f)^*$ of $f$, using (\ref{lcf}):
\begin{equation}
\label{esfer}
 f_{Sp}^*(x)=\int_0^{\infty} \chi_{\{f>s\}^*}(x)\ ds,
\end{equation}

\begin{equation}
\label{steiner} ({\cal S}_k f)^*(x)=\int_0^{\infty} \chi_{{\cal
S}_k(\{f>s\})}(x)\ ds.
\end{equation}

First, we observe that, by an easy change of variables in
(\ref{esfer}), we obtain that if $f^*$ denotes the classical
decreasing rearrangement of $f$,
\begin{equation}
\label{radial}
 f_{Sp}^*(x)=f^*(\sigma_n |x|^n), \ x\in \R^n.
\end{equation}
 In particular, the spherical
rearrangement  $f_{Sp}^*$, of a measurable function $f$ in $\R^n$,
is a radial decreasing function.

By a change of variables into spherical coordinates in $\R^k$, we
can write the last $k-$coordinates of $x\in \R^n$, $(x_{n-k+1},
\cdots, x_n)$ as $\rho\theta_{k-1}$, with $\rho>0$ and
$\theta_{k-1}\in \Sigma_{k-1}$ (the unit sphere in $\R^{k}$). Thus, using (\ref{radial}), we have that
\begin{equation}
\label{radial2}
 ({\cal S}_k
f)^*(x)=(f_{\bar{x}})^*(\sigma_k \rho^k),
\end{equation}
 where
$(f_{\bar{x}})^*$  is the classical decreasing rearrangement of
the function defined on $\R^k$ as follows: 
$f_{\bar{x}}(\cdot):=f(\bar{x},\cdot)$, with respect to  the
coordinates $(x_{n-k+1}, \cdots, x_n)\in \R^k$.

 Taking into account these considerations,
given a weight $v$ defined on $\R^n$, we can establish the
following formula for the functional defining the Lorentz spaces
$\Lambda^p_{{\cal S}_k}(v)$, with respect to the rearrangement given
by the Steiner symmetrization of order $k$.

\begin{proposition}
\label{cooresf}
 Let $0<p<\infty$. Given a weight $v$ defined on $\R^n$,
there exists another weight $\bar{v}$ on $\R^{n-k}\times \R^+$ such
that, for all measurable functions  $f$ on $\R^n$,
$$ \|f \|^p_{\Lambda^p_{{\cal S}_k}(v)}=\frac{1}{k\sigma_k}
\int_{\R^{n-k}} \int_0^{\infty} (f_{\bar{x}})^{*p}(s) \
\bar{v}(\bar{x},s) \ ds\ d\bar{x}.$$
\end{proposition}

\noindent{\bf Proof:} This is just a consequence of (\ref{radial2})
after a change into spherical coordinates $\rho\theta_{k-1}$, with
$\rho>0$ and $\theta_{k-1}\in \Sigma_{k-1}$. Then calling
$s=\sigma_k \rho^k$, the weight $\bar{v}$ associated to $v$ is
given by
\begin{equation}
\label{associate}
 \bar{v}(\bar{x},s):=\int_{\Sigma_{k-1}}
v(\bar{x},(s/\sigma_k)^{1/k} \theta_{k-1}) \ d\theta_{k-1}, \ \
(\bar{x},s)\in \R^{n-k}\times \R^+ .
\end{equation}
 $\hfill\Box$ \bigbreak

\begin{remark}{\rm
We remark that in the case of Steiner symmetrization of order
$k=1$ (the corresponding to one dimensional cross
sections), the associated weight to $v$ is just
$\bar{v}(\bar{x},y)=v(\bar{x},y)+v(\bar{x},-y), \ (\bar{x},y)\in
\R^{n-1}\times \R^+$.}
\end{remark}

Looking at the formula (\ref{satur}), by means of a change into
spherical coordinates and the use of (\ref{radial}), we can deduce
that, if $u$ is a decreasing function in $\R^+$,  then  the
following saturation formula for the spherical rearrangement hold:
\begin{equation}
\label{satur4}
 \sup_{\sigma} \int_{\R^n} |f(x)| \ u(\sigma(x)) \ dx
=\int_{\R^n} f_{Sp}^*(x) u(\sigma_n |x|^n) \ dx,
\end{equation}
 where the supremum is taken
over all   measure preserving transformations
$\sigma:\R^n\longrightarrow \R^{+}$.

All these facts lead us to establish the following
characterization of the normability of Lorentz spaces with respect
the spherical rearrangement with essentially the same proof as in
the classical case (see \cite{Lo}).

\begin{theorem}
Let $v$ be a weight in $\R^n$ and  $p\geq 1$. The following facts
are equivalent:

\begin{enumerate}
\item[(a)] The functional $\|\cdot \|_{\Lambda^p_{Sp}(v)}$ is a norm.
\item[(b)] For every $A, B\in \R^n$, $V((A\cup
B)^*)+V((A\cap B)^*)\leq V(A^*)+V(B^*).$
\item[(c)] The   weight $\bar{v}(s):=\displaystyle{\int_{\Sigma_{n-1}}
v((s/\sigma_n)^{1/n} \theta_{n-1}) \ d\theta_{n-1}}, \ \ s\in
\R^+,
 $ is a decreasing function.
\item[(d)] For all measurable functions $f$ in $\R^n$, the equality
$$\sup_{h_{Sp}^*=\bar{v}} \int_{\R^n} |f(x) h(x)| \ dx=
\int_{\R^{n}} f_{Sp}^*(x)\bar{v}(\sigma_n |x|^n) \ dx $$ holds.
\end{enumerate}
\end{theorem}

\noindent{\bf Proof:}
 Theorem \ref{tconcavity} gives us that (a) implies
(b).

Assume that (b) holds, and consider $0<\varepsilon <a\leq b$ and the following 
sets in $\R^n$,
$$ A=B(0,a), \\ B=B(0,b)\setminus B(0,\varepsilon).$$
Then,
$$ A=A^*, \ B^*=B(0,(b^n-\varepsilon^n)^{1/n}), \ (A\cup B)^*=
B(0,b), \ (A\cap B)^*=B(0,(a^n-\varepsilon^n)^{1/n}).$$

 Condition (b) implies that
 $$V(B(0,b))-V(B(0,(b^n-\varepsilon^n)^{1/n}))\leq V(B(0,a))-V(B(0,(a^n-\varepsilon^n)^{1/n})).$$

 After a change into spherical coordinates $(\rho, \theta_{n-1})\in \R^+\times \Sigma_{n-1}$
 and calling $s=\sigma_n \rho^n$, we obtain that the above condition
 can be written as
 $$ \int_{\sigma_n(b^n-\varepsilon^n)}^{\sigma_n b^n} \bar{v}(s) \ ds
 \leq \int_{\sigma_n(a^n-\varepsilon^n)}^{\sigma_n a^n} \bar{v}(s) \
 ds. $$

Dividing both sides by $\sigma_n\varepsilon^n$ and letting $\varepsilon
\rightarrow 0$, we obtain $\bar{v}(\sigma_n b^n)\leq
\bar{v}(\sigma_n a^n)$; that is, $\bar{v}$ is a decreasing function
of $s$.

That condition (c) implies (d) is equality (\ref{satur4}). Finally, we observe that Theorem~\ref{saturation} proves that
condition (d) implies (a). $\hfill\Box$ \bigbreak

Similarly, in order to study when the functional
$\|\cdot\|_{\Lambda^p_{{\cal S}_k}(v)}$ is a norm, we observe
that, due to Proposition \ref{cooresf}, the condition is reduced
to the fact that the associated weight $\bar{v}(\bar{x},s)$,
defined in (\ref{associate}), must be a decreasing function in
$s$, and, also with essentially the same proof, we can establish
the following characterization.

\begin{theorem}
Let $v$ be a weight in $\R^n$,   $p\geq 1$ and $k\geq 1$   an
integer. The following facts are equivalent:
\begin{enumerate}
\item[(a)] The functional $\|\cdot \|_{\Lambda^p_{{\cal
S}_k}(v)}$ is a norm.

\item[(b)] For every $A, B\subset \R^n$, $V({\cal S}_k(A\cup
B))+V({\cal S}_k(A\cap B))\leq V({\cal S}_k(A))+V({\cal S}_k(B)).$

\item[(c)] The   weight $\bar{v}(\bar{x},s)$ defined in (\ref{associate}) is a
decreasing function in the variable $s$.

\item[(d)] For all measurable functions $f$ in $\R^n$, the equality
$$\sup_{({\cal S}_k h)^*=\bar{v}} \int_{\R^n} |f(x) h(x)| \ dx=
\int_{\R^{n-k}\times \R^+}({\cal S}_k f)^*(\bar{x},s) \
\bar{v}(\bar{x},s) \ d\bar{x} \ ds$$ holds.

\end{enumerate}
\end{theorem}

\address{

\noindent
Santiago Boza\\
Dept. Appl. Math. IV\\
EPSEVG\\ Polytechnical University of Catalonia\\ E-08880 Vilanova
i Geltr\'u, SPAIN\ \ \ {\sl E-mail:}
 {\tt boza@ma4.upc.edu}

\medskip

\noindent Javier Soria\\ Dept. Appl. Math. and Analysis
\\ University of Barcelona\\ E-08007 Barcelona,
 SPAIN\ \ \ {\sl E-mail:}
 {\tt soria@ub.edu}
}
\end{document}